\input amstex
\documentstyle{amsppt}
\magnification=1200
\pageheight{7.5in}
\pagewidth{5.4in}
%\hcorrection{.25in}
%\advance\vsize-.75in
\input epsf
\NoBlackBoxes
\topmatter
\title
Positive Diagrams for Seifert Fibered Spaces
\endtitle
\author
John Hempel
\endauthor
\affil
Rice University, Houston, Texas
\endaffil
\address
Rice University, Houston, Texas
\endaddress
\email
hempel\@math.rice.edu
\endemail
\endtopmatter

\document 

%{\bf Preliminary version 6/99}

\smallpagebreak

{\bf 0. Introduction.} By a {\it   (Heegaard) diagram} we mean a triple $(S; X, Y)$ where $S$
is a  closed,  connected, orientable surface and $X$ and $Y$ are
compact  1-manifolds in $ S$.

 A  diagram gives rise to a 3-manifold
$$
M= S \times [-1,1]\bigcup_{X \times -1} \text{2-handles} \bigcup_{Y \times 1}
 \text{2-handles} \bigcup \text{3-handles}
$$
 obtained by adding 2-handles to $S \times [-1,1]$ along the curves of
 $X \times -1$ and $Y \times 1$ and then adding 3-handles
 along all resulting 2-sphere boundary components. The decomposition
 of $M$ by $S \times 0$ is the associated {\it (Heegaard)  splitting}
 of $M$ and the genus of $S$ is called the {\it genus} of the splitting.

 A {\it positive diagram} is a diagram in which $S$, $X$, and $Y$ are
 oriented  and the
 intersection number $<X,Y>_p$ of $X$ with $Y$ is $+1$ at each point
 $p \in X\cap Y$.  
 Every compact, oriented  3-manifold (with no 2-sphere boundary
 components) is represented by a positive diagram. One can start with
 an arbitrary diagram for the manifold and eliminate negative
 crossings by adding trivial handles. The new curves associated with
 each  trivial handle can be oriented so as to introduce only positive
 crossings. In the process one negative crossing is replaced by three
 positive crossings and the genus of the associated  splitting
is increased by one.

Positive diagrams were introduced by Montesinos [M] who observed that
they have a nice encoding by the permutations on the $d = \#(X \cap
Y)$ intersection points given by flowing along $X$ (respectively $Y$)
from one point to the next, and noted that these permutations are also
the monodromy representation for the manifold as a branched cover of
$S^3$ branched over a fixed (universal) branch set. See [H3] for  more
information.

It is natural to wonder about the relation between the {\it Heegaard
genus},  $hg(M)$, and the {\it positive Heegaard genus},  $phg(M)$, of $M$;
where these are defined respectively by minimizing the genus among all
splittings of $M$ and among all splittings determined by a positive
diagram.
From the above it might seem that the difference $phg(M)-hg(M)$ could
be arbitrarily large; as one might have to increase the genus by as
much as $\#(X \cap Y)/2$ to change all  the negative crossings.  However this cannot be supported by group theory alone
in the following sense. Assuming that at least one of the sides of the
splitting is a handlebody, we get a presentation for $\pi_1(M)$ with
generators dual to the attaching curves for the handlebody and
relators comming from the attaching curves for the other side. This will be a {\it positive presentation}: the relators are words in
positive powers of the generators.

\proclaim{Fact} If a group $G$ has a presentation with $n$ generators
and $m$ relators, then it has a positive presentation with $n+1$
generators and $m+1$ relators.
\endproclaim
\demo{Proof} Add a generator $x_{n+1}$ and a relator $x_1x_2 \dots
x_nx_{n+1}$.
In every other relator, and for each $i=1,2,\dots,n$ replace each
occurrance of $x_i^{-1}$ by $x_{i+1} \dots x_nx_1 \dots x_{i-1}$. \qed
\enddemo

If one starts with a presentation comming from a Heegaard splitting,
the  positive presentation obtained as above will not, in general, correspond
to a splitting -- one cannot represent the new/modified relations by
disjoint, embedded simple closed curves in the boundary of a
handlebody. 

To my knowledge the only class of  3-manifolds (including Heegaard genus
$>2$) for which there is exact determination of Heegaard genus
[BZ],[MS] are the Seifert manifolds (see Section 1 for the notation we
use for them). For these we prove:

\proclaim{Theorem A} Let $M$ be a  closed, orientable Seifert fibered
space with orientable base space. Then $phg(M) = hg(M)$ except for the
following cases:

(1) $g=0$; $m \ge 4$ is even; Seifert invariants $ = 1/2,1/2, \dots
    ,1/2, n/2n + 1\quad (n \ge1)$; $e=m/2$. In these cases
    $phg(M)\le hg(M)+1=m-1$, and $phg(M)= hg(M)+1$ if $m=4$ or $n >1$.

(2) $g>0$; and ( $m=0$;  $e=\pm 1$) or ($m=1$, Seifert invariant $=
    1/\alpha$, $e=0$) or ($m=1$, Seifert invariant $=(\alpha
    -1)/\alpha$, $e=1$). In these cases $2g+1 = hg(M)+1 \le phg(M) \le
    hg(M)+2$.

(3) $g>0$; $m \le 2$; and not in case (2) above. Then $
    min\{2g+1,2g+m-1\}=hg(M) \le phg(M)\le hg(M)+1$.
\endproclaim

In cases (1) and (2) of Theorem A the minimal genus splittings are
``horizontal splittings''. Not many Seifert manifolds have horizontal
splittings and they are rarely of minimal genus. The next theorem
gives the remaining exceptions.

\proclaim{Theorem B} Let $M$ be a closed, oriented, Seifert fibered
space with orientable base space which is not included in Theorem A
and which has a horizontal splitting realizing $hg(M)$. Then $M$ has
$g=0$, $m=3$, $e=1$, and Seifert invariants $(1/2,1/3,n/(6n \pm 1)$,
or $(1/2,1/4,n/(4n \pm 1))$, or $(1/3,1/3,n/(3n \pm 1)$ for some $n
\ge 1$.

Except for the cases $(1/2,1/3,1/5), (1/2,1/4,1/3), (1/3,1/3,1/2),$
and $(1/2,1/3,1/7)$ the horizontal splitting is not represented by a
positive diagram.
\endproclaim

In the first three exceptions of Theorem B the splitting is also a
vertical splitting and so is represented by a positive diagram by
Theorem C, below. We are unsure about the fourth case.

\proclaim{Theorem C} Every vertical splitting of an orientable Seifert
fibered space with base space $S^2$ is represented by a positive
diagram.
\endproclaim

The proofs of these theorems depend on the result [MS] that splittings
of Seifert manifolds are either ``vertical'' or ``horizontal'' and
are organized as follows. In section 2 we discuss vertical splittings
and prove Theorem C. In section 3 we apply Theorem C and some lifting
arguements to establish the upper bounds of Theorem A (Corollary 3.4)
and the equality $phg(M) = hg(M)$ when $M$ has at least three singular
fibers and the base surface is not $S^2$ (Theorem 3.6).
In section 4 we discuss horizontal splittings and give (Theorem 4.1)
the Seifert manifolds for which these horizontal splittings can be of
minimal genus.  We prove (Theorem 4.2) that, with the indicated
exceptions, these minimal genus horizontal splittings are not
represented by positive diagrams. This gives the lower bounds of
Theorem A and the proof of Theorem B. In section 5 we discuss some of
the questions left open.

{\bf 1. Preliminaries.} A {\it compression body} is a space built as
follows. Take a closed, orientable surface $S$, attach 2-handles to $S
\times [0,1]$ along curves in $S \times 1$, and fill in any resulting
2-sphere boundary components with 3-balls. The image of $S \times 0$
is called the {\it outer boundary component} of the compression body. 

A {\it (Heegaard) splitting} of a compact,
orientable 3-manifold $M$ is a representation of $M$ as the union of
two compression bodies (handlebodies if $M$ is closed) meeting in a
common outer boundary component. Formally, it is a pair $(M,S)$ where
$S \subset M$ is a closed, orientable surface which separates $M$ so
that the closure of each component of $M-S$ is a compression body with
$S$ as outer boundary component.

A {\it (Heegaard) diagram} is a description of a Heegaard splitting by
designating the common outer boundary component $S$ and sets of
attaching curves in $S$ for the 2-handles  each of the compression
bodies  is to bound. Formally it is a triad $(S;X,Y)$ where $S$ is a closed,
orientable surface and $X$ and $Y$ are compact 1-manifolds in $S$

Often we will prefer to work with {\it oriented} objects. So an
oriented Heegaard splitting is a pair $(M,S)$ of oriented manifolds
(as above) and equivalence of such will be an  orientation preserving
homeomorphism of pairs. An oriented diagram will be an oriented triad
$(S;X,Y)$, and will determine an oriented splitting by the convention
that the positive normal to $S$ in $M$ points toward the $Y$-side of
the splitting.

We adopt the following notation and conventions for a closed, oriented
Seifert fibered space $M$ :

$\bullet \quad g \ge 0$ will denote the genus of the base surface.

$\bullet \quad m \ge 0$ will denote the number of singular fibers.

$\bullet \quad e \in {\Bbb Z} $ will denote the Euler number.

$\bullet \quad (\alpha_i,\beta_i) \in {\Bbb Z} \times {\Bbb Z} ; 1 \le i \le m$ will
denote the Seifert invariants.

It is understood that $\alpha_i$ and $ \beta_i$ are relatively prime
with $0<\beta_i<\alpha_i$. Some times we will write these as fractions
$\beta_i/\alpha_i$.

We refer  the data as above the {\it normalized} invariants for
$M$. They are unique to $M$ which is obtained by oriented Dehn filling
on a product $F \times S^1$ according to the formula
$$
x_i^{\alpha_i}t^{\beta_i} = 1; i=1, \dots, m
$$
$$
[a_1,b_1] \dots [a_g,b_g]x_1 \dots x_mt^e = 1
$$
where $F$ is an oriented surface with genus $g$ and $m+1$ boundary
components.

We sometimes find it more convenient to work with a non-unique form of
the invariants. In particular, finding a positive diagram depends on
making a suitable choice of these invariants. So if we take an
oriented surface $F$ of genus $g$ and with $r \ge 1$ oriented boundary
components $x_1, \dots, x_r$ and we do Dehn filling on $F \times S^1$
according to the formula
$$
x_i^{\alpha_i}t^{\beta_i'} =1; i = 1, \dots, r
$$
where $\alpha_i \ge 1$ and $\beta_i'$ is prime to $\alpha_i$,
the resulting manifold is said to have {\it non-normalized} invariants
$\{g; \beta_1'/\alpha_1, \dots, \beta_r'/\alpha_r\}$.

One can change these invariants by doing twists on vertical annuli in
$F \times S^1$. This way one can show

\proclaim{1.1 Proposition} The normalized invariants for the Seifert
fibered space with non-normalized invariants $\{g; \beta_1'/\alpha_1,
\dots, \beta_r'/\alpha_r\}$ are

genus = $g$

$m=\#\{i:\alpha_i > 1\}$

Seifert invariants $\beta_i/\alpha_i$ when $\alpha_i >1$ and
$\beta_i$ is the least positive residue of $\beta_i'$ modulo
$\alpha_i$.

$e= -\sum[\beta_i'/\alpha_i]$
\endproclaim

{\bf 2. Vertical splittings.} It is known [MS] that every splitting of
an orientable Seifert manifold with orientable base space is either
{\it horizontal} or {\it vertical} . All Seifert manifolds have
vertical splittings, but most do not admit horizontal
splittings. Theorem 0.3 of [MS] describes, in terms of the Seifert
invariants, those  Seifert manifold which have  horizontal
splittings. We review these definitions here.

Let $M$ be a Seifert manifold  with base surface $B$ and projection
$f:M \to B$. Suppose we have a cell decomposition of $B$ such that $B
= D \cup E \cup F$ where each of $D,E,F$ is a disjoint union of closed
2-cells of the decomposition, each component of $D$ and of $E$
contains at most one singular point, which is an interior point,  each
component of $F$ is a square containing no singular point and having
one pair of opposite sides in $D$ and the other pair in $E$,  $Int(D)
\cap Int(E) = Int(D) \cap Int(F) = Int(E) \cap Int(F) = \emptyset$,
and $D \cup F$ and $E \cup F$ are connected. See Figure 1.

\bigpagebreak
\epsfysize=1.15in
\centerline{\epsffile{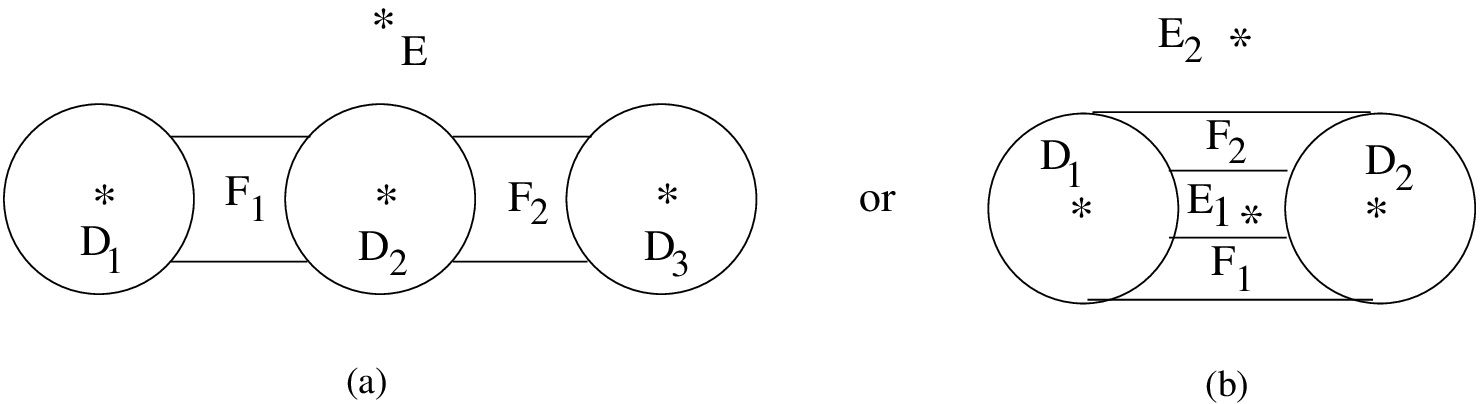}}
\botcaption{Figure 1}\endcaption
\bigpagebreak

Then $f^{-1}(D)$ is homeomorphic to $D \times S^1$,  the same holds
for $E$ and $F$, and we consistently fix such identifications. Let $V_1 = D \times
S^1\cup F \times [0,1/2]$ and $V_2 = E \times S^1 \cup F \times [1/2,1]$; where $S^1 = [0,1]/0 \sim 1$.
Put $S=V_1 \cap V_2 =\partial(V_1) = \partial(V_2)$. Then $(M,S)$ is a splitting of $M$ of genus
$$
g = \beta_0(D)+\beta_1(D \cup F)= \beta_0(E) + \beta_1(E \cup F) = 1 + \beta_0(F)
$$
$$
=2g(S)+\beta_0(D)+\beta_0(E)-1
$$
 which we call a {\it vertical splitting}.  These calculations
 immediately give

\proclaim{2.1 Proposition} A vertical splitting of a closed,
 orientable Seifert manifold over an orientable surface of genus $g$
 and with $m$ singular fibers has genus at least
$$
max\{2g+1,2g+m-1\}
$$
\endproclaim

{\bf Proof of Theorem C.} We may assume $m \ge 3$. We take the vertical splitting as described
above  where $D$ has components $D_1, \dots, D_{m_1}$ and $E$ has
components $E_1, \dots, E_{m_2}$; where $m_1+m_2=m$. Let $\mu_i$ and  $\nu_j$ be the
positively oriented boundaries of $D_i$ and $E_j$ respectively. $M$ is
obtained from $(\partial D \cup \partial E \cup F)\times S^1$ by
attaching solid tori according to formula
$$
\mu_i^{\alpha_i}t_i^{\beta_i}=1;\quad i=1,\dots,m_1
$$
$$
\nu_j^{\alpha_{m_1+j}}s_j^{\beta_{m_1+j}}=1;\quad j=1,\dots,m_2
$$
corresponding to some non-normalized coordinates $\{\beta_i/\alpha_i\}$. The curves $t_i$ ($s_j$)
are positively oriented vertical curves in $f^{-1}(D_i \cap E)$ ($f^{-1}(E_j \cap D)$). Recall that, by 2.1, we
may change the $\beta_i$ modulo $\alpha_i$ as long as we don't change
$\sum[\beta_i/\alpha_i]$. We will impose more conditions on these
invariants later.

Let $\Gamma_D$ (respectively $\Gamma_E$) be the graph whose vertices
correspond to the components of $D$ (respectively of $E$)
and whose edges correspond to the components of $F$. There are
natural embeddings $\Gamma _D \subset D \cup F, \Gamma_E \subset E
\cup F$.

The meridian disks  for $V_1$ will be:

(i) meridian disks for the filling solid $f^{-1}(D_i)$; together with

(ii) vertical disks $A_p \subset F_p \times [0,1/2]$ corresponding to
those components $F_p$ of $F$ not in a fixed maximal tree of
$\Gamma_D$. These will be of the form $c_p \times [0,1/2]$, where
$c_p$ is an arc in $F_p$ separating its edges lying in $D$.

Similarly the meridian disks for $V_2$ will be meridian disks for the
$f^{-1}(E_j)$ and vertical disks $B_q \subset F_q \times [1/2,1]$
corresponding to components $F_q$ of $F$ not in a maximal tree in
$\Gamma_E$.

We need the following conditions.

(1) Each $A_p$ has vertical sides in some $f^{-1}(E_{j_1})$ and some
    $f^{-1}(E_{j_2})$. We require that $\beta_{m_1+j_1}$ and
    $\beta_{m_1+j_2}$ have opposite sign. Same for the $\beta$'s on opposite sides of each $B_q$.

(2) $A_p \cap B_q = \emptyset$ for all $p,q$.

(3) $s_j$ is to be chosen in $f^{-1}(D_i \cap E_j)$ with $\beta_i >
    0$. Similarly $t_i$ is to be chosen in $f^{-1}(D_i \cap E_j)$ with
    $\beta_{m_1+j} > 0$

First we show that these conditions will produce a positive
diagram. Their justification will  be given afterwards. Figure 2
illustrates this construction using the decomposition of Figure 1b.

\bigpagebreak
\epsfysize=2.6in
\centerline{\epsffile{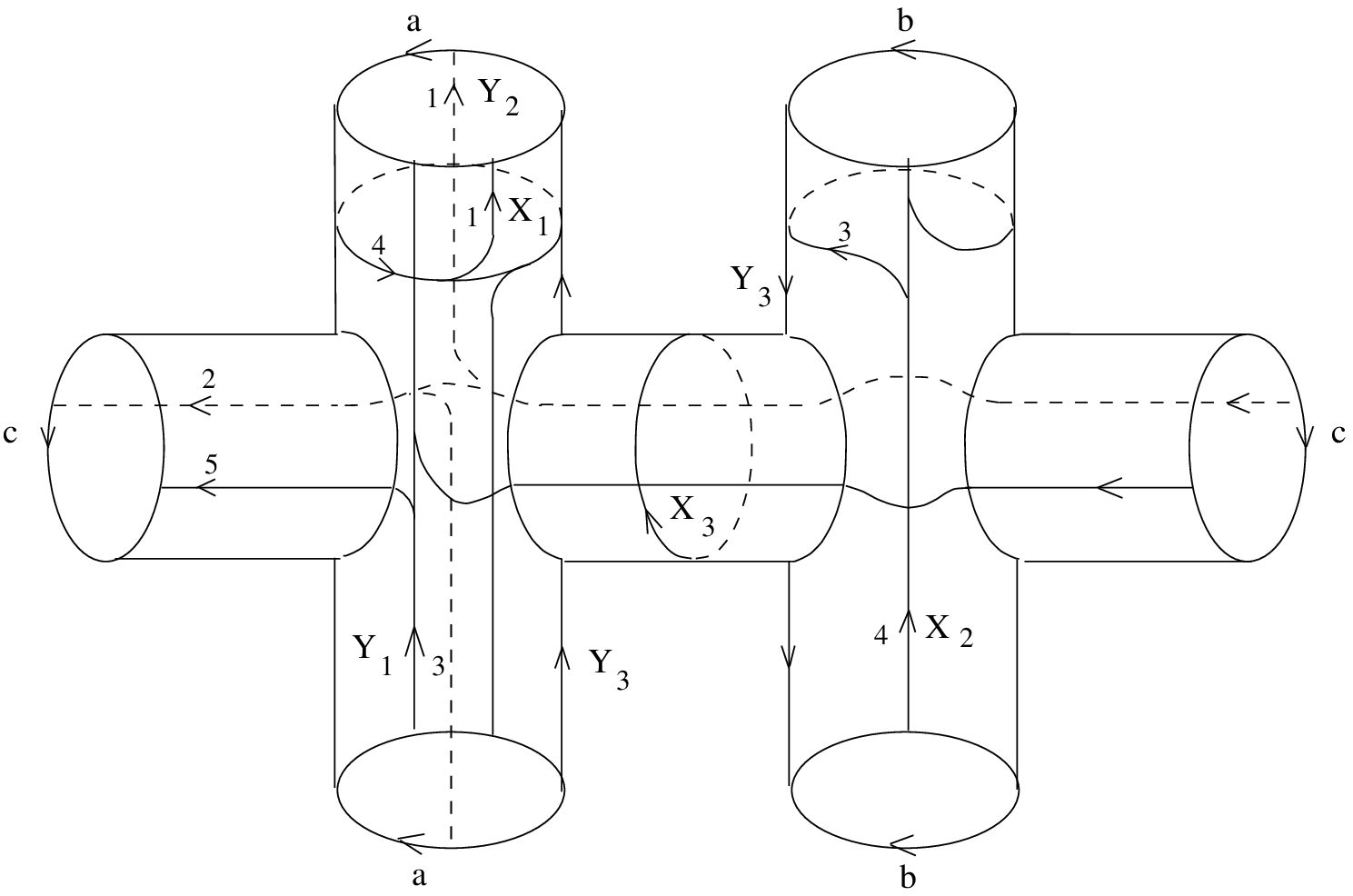}}
\botcaption{Figure 2}
 Positive diagram for Seifert manifold
 $(g=0;1/4,2/3,3/5,1/2; e=5)$; for which we use non-normalized
 invariants $(1/4,-4/3,3/5,-5/2)$. Curves are represented as weighted
 train tracks.
\endcaption
\bigpagebreak

The meridian curve $X_i$ for $f^{-1}(D_i)$ is chosen in a neighborhood
of $t_i \cup \mu_i \times 3/4$ and oriented so that it is homotopic in
$f^{-1}(\partial D_i)$ to $\mu_i^{\alpha_i}t_i^{\beta_i} $ or
$\mu_i^{-\alpha_i}t_i^{-\beta_i}$ according as $\beta_i >0$ or $\beta_i
<0$. We similarly choose the meridian $Y_j$ of $f^{-1}(E_j)$.

Now we may suppose $S = \partial V_1 = \partial V_2$ is oriented so
that
$<\mu_i, t> = +1$ for every positively oriented vertical curve $t$ in
$f^{-1}(\partial D_i)$. It follows that $<t,\nu_j> = +1$ for every
positively oriented vertical curve $t$ in
$f^{-1}(\partial E_j)$. Now $X_i \cap Y_j$ lies in a neighborhood of
$((\mu_i \times 3/4) \cap s_j) \cup (t_i \cap (\nu_j \times
1/4))$. The first (second) of these terms is nonempty only if $s_j
(t_i) \subset  f^{-1}(D_i\cap
E_j)$. Condition (3) then assures us that $<X_i,Y_j> = +1$ at each
crossing point.

By (2) a vertical $A_p$ will only meet the two $Y$ curves $Y_{j_1}$ and
$Y_{j_2}$. By condition (1) $\partial A_p$ can be oriented to meet both
of these positively at each point. Similarly each $\partial B_q$ can
be oriented so that $X$ crosses it positively at every point.

To justify condition (1), let T be the subgraph of $\Gamma_E$ dual to
the edges of $\Gamma_D$ not in a maximal tree $T_D$ in $\Gamma_D$. $T$
is connected as $S^2-T_D$  is connected.  $T$ must
be a tree; otherwise some loop in $T$ meets a loop in $\Gamma_D$ in a
single point. This is impossible in $S^2$. Thus we can put the
vertices of $T$ into two classes according to whether the simplicial
distance in $T$ to a fixed vertex is even or odd. Each edge of $T$
joins a vertex of one class to one of the other class. We can then
choose the $\beta's$ corresponding to one class positive and to the
other negative. This is possible as we can change the $\beta$'s
(modulo the corresponding $\alpha$'s) as long as we keep $\sum
[\beta_i/\alpha_i]$ constant.

Note that the above choices of $T_D \subset \Gamma_D$ and $T \subset
\Gamma_E$ (which is clearly maximal) satisfy (2); since $(\Gamma_D -
T_D)\cap (\Gamma_E-T) = \emptyset$

To get (3) it suffices to show that each $E_j$ meets some $D_i$ with
$\beta_i >0$ (and vice versa). Note that $E_j$ 
is contained in a unique component $E_j'$ of  $S^2 - \Gamma_D$. If $E_j$ only met
$D_i$'s with $\beta_i < 0$, then, by (1), all the edges of $\Gamma_D$
in $\partial E_j'$ lie in the maximal tree. This is impossible unless
$\Gamma_D$ is itself a tree (and $E$ has a single component). This
last possibility is easily handled.

{\bf 3. Lifting to covers.} The following result comes from lifting a
positive diagram for the base.

\proclaim{3.1 Lemma} If $M$ is represented by a positive diagram of
genus $g$ and $p:\tilde M \to M$ is a $\lambda$-sheeted covering
space, then $\tilde M$ is represented by a positive diagram of genus
$\tilde g = \lambda (g-1)+1$.
\endproclaim

The next result describes how to lift a Seifert fibration to a
particular kind of cover.

\proclaim{3.2 Lemma} Let $M$ be an orientable Seifert fibered space
over an orientable surface of genus $g$ and with non-normalized
invariants $\beta_1/\alpha_1, \dots, \beta_r/\alpha_r; r >0$. So $M$
is obtained by Dehn filling on $F \times S^1$ where $F$ is an oriented
surface of genus $g$ and with oriented boundary components $x_1,
\dots, x_r$. Let $p:\tilde F \to F$ be a $\lambda$-sheeted covering so
that for each $i = 1, \dots,r$ $p^{-1}(x_i)$ has components $\tilde
x_{i,1}, \dots, \tilde x_{i,r_i}$ with $p|\tilde x_{i,j}:\tilde
x_{i,j} \to x_i$ a $b_{i,j}$-sheeted cover and with $b_{i,j} $
dividing $\alpha_i$. Put $\alpha_{i,j} = \alpha_i/b_{i,j}$.

Then $p \times 1:\tilde F \times S^1 \to F \times S^1$ extends to a
$\lambda$-sheeted covering space $\tilde M \to M$ where $\tilde M$ is
the orientable Seifert fibered space over the orientable surface of
genus
$$
\tilde g = \lambda(g-1)+1+(r\lambda-\sum r_i)/2
=\lambda(g-1)+1+\sum (b_{i,j}-1)/2  
$$
 with non-normalized invariants
$\{\beta_i/\alpha_{i,j}; \quad i= 1, \dots, r, j=1,\dots,r_i\}$.

\endproclaim

\proclaim{3.3 Theorem} Let $\tilde M$ be an orientable Seifert fibered
space over an orientable surface of genus $g$ with at most  three
singular fibers. Then $\tilde M $ is a $(2g+1)$-sheeted cover of an
orientable Seifert space over $S^2$ with three singular fibers.
\endproclaim

\demo{Proof} We represent $\tilde M$ with non-normalized invariants
$\beta_i/\alpha_i; \quad i=1,2,3$ (some $\alpha_i=1$ if $m < 3$).
 By Lemma 3.5
below, with $\lambda = 2g+1$, we can find integers $\beta_i^*$ with
$\beta_i^* \equiv \beta_i \text{ mod } \alpha_i$, $(\lambda,\beta_i^*)=1$, and
$[\beta_1^*/\alpha_1] + [\beta_2^*/\alpha_2]+[\beta_3^*/\alpha_3]
=[\beta_1/\alpha_1] + [\beta_2/\alpha_2]+[\beta_3/\alpha_3]$.

Let $F$ be a sphere with three holes. Since $\lambda$ is odd, there is
a $\lambda$-sheeted cyclic
covering space $p:\tilde F \to F$ such that the inverse image of each
boundary component of $F$ is connected.

Now apply Lemma 3.2 to the Dehn filling of $F \times S^1$ with
non-normalized invariants $\beta_i^*/\lambda \alpha_i;\quad  i=1,2,3$. The
$\lambda$-sheeted cover of $M$ given by 3.2 is our given $\tilde M$. \qed
\enddemo

By Theorem C an orientable Siefert fibered space over $S^2$ with three
singilar fibers is represented by a positive diagram of genus
two. Applying Lemma 3.1 and Theorem 3.3 to this gives

\proclaim{3.4 Corollary} An orientable Seifert fibered space over an
orientable surface of genus $g$ with at most three singular fibers
is represented by a positive diagram of genus $2g+2$.
\endproclaim

\proclaim{3.5 Lemma}Given $(\alpha_i, \beta_i)\in {\Bbb Z} \times
{\Bbb Z}$ with $\alpha_i \ge 1$ and $(\alpha_i, \beta_i)=1;\quad
i=1,\dots,n$, and given odd $\lambda \in {\Bbb Z}$ there exist
$\beta_i^* \in {\Bbb Z};\quad i = 1,\dots,n$ satisfying

(1) $\beta_i^* \equiv \beta_i \text{ mod } \alpha_i$

(2) $(\beta_i^*,\lambda)=1$

(3) $\sum[\beta_i^*/\alpha_i] = \sum[\beta_i/\alpha_i]$
\endproclaim

\demo{Proof} We induct on the number of distinct prime factors of $\lambda$.

Let $\lambda = p^r$ where $p$ is an odd prime. If $(p,\beta_i)=1$ for
all $i$, we are already done. So suppose $p|\beta_i$ if and only if $1
\le i \le m$. Then $(p,\alpha_i)=1 ; \quad 1 \le i \le m$; so
$(p,\beta_i+k\alpha_i)=1; \quad 1 \le i \le m$ if $1 \le |k|<p$.

If $m=1$, then $p$ cannot divide both  $\beta_2 +\alpha_2$ and
$\beta_2 - \alpha_2$; otherwise $p$ divides $2\beta_2$ ; hence
$p$ divides $\beta_2$. So, say, $(p, \beta_2-\alpha_2)=1$. Put
$\beta_1^*=\beta_1+\alpha_1, \quad \beta_2^* = \beta_2 - \alpha_2$, and
$\beta_i^*=\beta_i; \quad i\ge 3$. 

If $m$ is even, put $\beta_i^*=\beta_i+\alpha _i; \quad 1\le i \le m/2$,
$\beta_i* = \beta_i - \alpha_i ; \quad m/2 < i \le m$, and $\beta_i^* =
\beta_i; \quad m < i$.

If $m > 1$ is odd, put $\beta_1^* = \beta_1+2\alpha_1$,\quad
$\beta_i^*=\beta_i+\alpha _i; \quad 2\le i \le (m-1)/2$, \quad $\beta_i^* = \beta_i
- \alpha_i ; \quad (m-1)/2 < i \le m$, and $\beta_i^* =\beta_i; \quad m < i$.

In each case $\{\beta_i^*\}$ satisfy (1), (2), and (3).

Now suppose $\lambda = \lambda'\lambda''$ where $(\lambda',
\lambda'')=1$, and that we have $\{\beta_i'\}$ and $\{\beta_i''\}$ satisfying
(1), (2), and (3) relative to $\lambda'$ and $\lambda''$ respectively.

Since $\alpha_i = (\alpha_i \lambda', \alpha_i \lambda'')$ divides
$\beta_i' - \beta_i'' $ we can use  the Chinese remainder theorem to
find $\beta_i^*$ satisfying
$\beta_i^* \equiv \beta_i'$ mod $\alpha_i \lambda'$ and $\beta_i^* \equiv
\beta_i''$ mod $\alpha_i \lambda''$.

 Then the $\{\beta_i^*\}$ satisfy (1), (2), and
(3$'$)  $\sum[\beta_i^*/\alpha_i] =  \sum[\beta_i/\alpha_i] +k\lambda $,
for some $k$.  Put $\beta_1^{**} =
\beta_1^*-k\alpha_1\lambda$ and $\beta_i^{**} = \beta_i^*; \quad i\ge
2$. The $\{\beta_i^{**}\}$ satisfy (1), (2), and (3). \qed
\enddemo

\proclaim{3.6 Theorem} Let $M$ be a closed oriented Seifert fibered
space over an orientable surface of genus $g >0$ and with $m \ge 3$
singular fibers. Then $phg(M)=hg(M) = 2g+m-1$.
\endproclaim

\demo{Proof} The fact that $hg(M) = 2g+m-1$ is established in
[BZ]. Note that with these assumptions the minimal genus splitting is
always a vertical splitting.  The case $m=3$ follows directly from
Corollary 3.4.

So we assume that $m >3$ and that $M$ is defined by non-normalized
invariants $\beta_i/\alpha_i;\quad i = 1, \dots,m$ with $\beta_i <0$ for
$i>3$.

Let $M'$ be the Seifert fibered space over the surface of genus $g$
with three singular fibers determined by the non-normalized invariants
$\beta_i/\alpha_i; i = 1,2,3$.

By 3.4 $M'$ is represented by a positive diagram $(S;X,Y)$ of genus
$2g+2$. From the proof of 3.4 (and its predicessors) there is a
positively oriented regular fiber $t \subset S$ so that $<X,t>_p =+1$
($<Y,t>_q = -1$) for each point $p \in X\cap t$ ($q \in Y \cap t$),
and so that there is an interval $J \subset S$ with $t\cap Y \subset
J$ and $J\cap X = \emptyset $ 

Choose parallel copies $t_4, \dots, t_m$ in $S$ which meet $X$ and $Y$
as $t$ does and push these to the
$Y$ side, $V_2$, of the splitting by disjoint annuli $t_i \times [0,1]$, where
$t_i = t_i \times 0$.

Now $M$ is obtained from $M'$ by surgery on the curves $t_i \times 1;\quad
i \ge 4$; so  ``small'' regular neighborhoods of these curves are
replaced by solid tori $W_i;\quad i \ge 4$. We tube these to $V_1$ along
regular neighborhoods of the $J_i \times [0,1]$ to obtain a handlebody
$V_1^*$ which will be half of the splitting for $M$. The other half,
$V_2$, is homeomorphic to the result of removing an open regular
neighborhood of $\bigcup_{i \ge 4}(J_i \times [0,1] \cup t_i \times 1)$.

We get a diagram for this splitting as follows. We add to $X$ meridian
curves for the $W_i$. These are oriented as $
\mu_i^{-\alpha_i}t_i^{-\beta_i}$ for a positively oriented regular
fiber $t_i$ and transversal $\mu_i$ in $\partial W_i$.

We modify the curves of $Y$ near where they cross the $J_i$ so as to
run around the other edge of a regular neighborhood of $J_i \times
[0,1]$. See Figure 3. To these we add the curves $((t_i - J_i)\times
[0,1])\cap \partial V_2^*; \quad i \ge 4$.

\bigpagebreak
\epsfysize=3.2in
\centerline{\epsffile{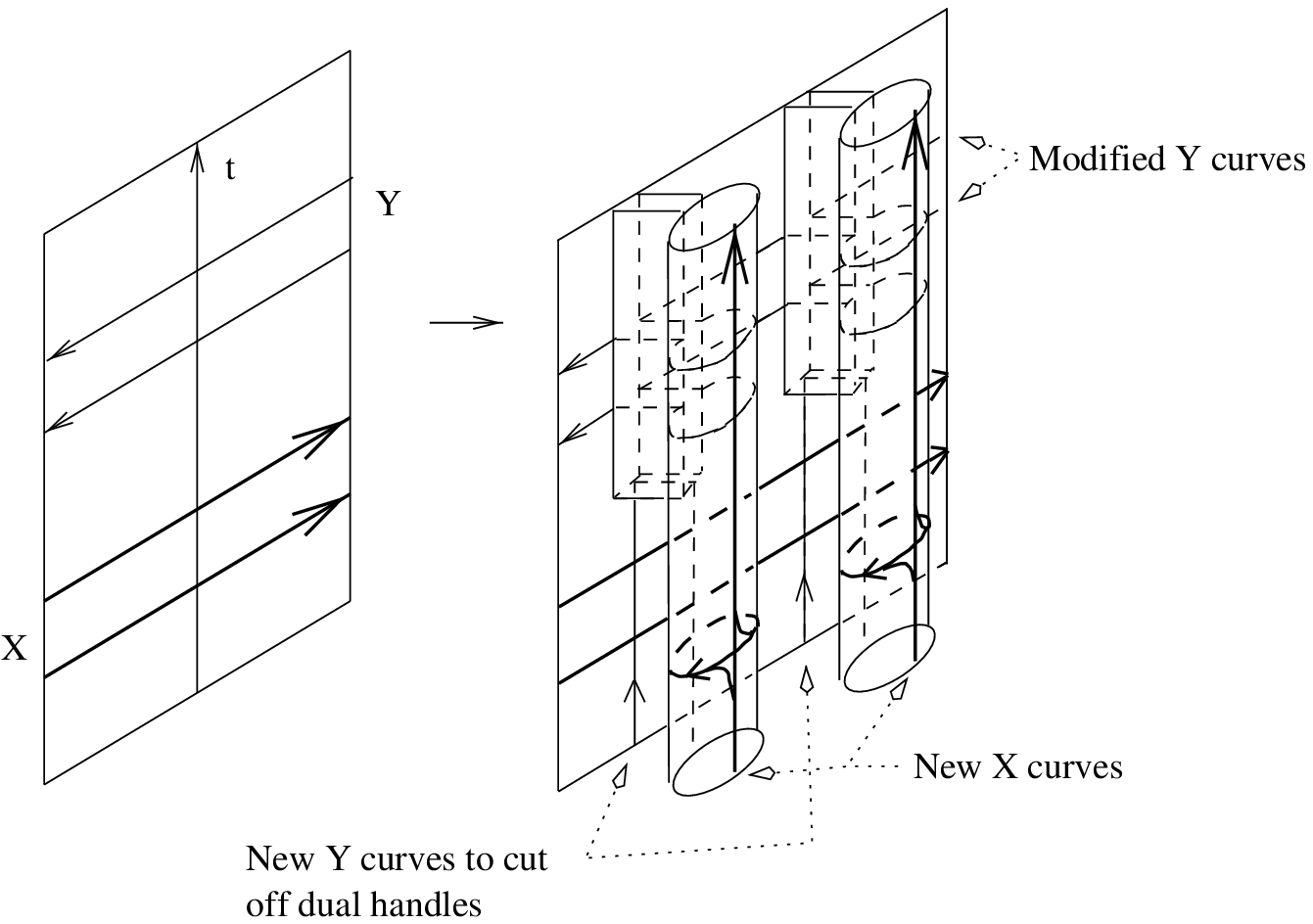}}
\botcaption{Figure 3}
 Modifying a positive diagram to accommodate more singular fibers.         
\endcaption
\bigpagebreak

One verifies that this construction can be made so that  all new intersections will be positive. \qed
\enddemo

{\bf 4. Horizontal splittings. } Let $F \ne B^2$ be a compact, connected,
orientable surface with one boundary component. Consider a surface
bundle over $S^1$
$$
N= F \times [0,1]/(x,0) \sim (\phi(x),1)
$$
where $\phi: F \to F$ is an orientation preserving homeomorphism with
$\phi | \partial(F) = 1$. Note that $\lambda = \partial(F) \times 0$
and $\mu = x_0\times [0,1]/\sim \quad (x_0 \in \partial(F))$ form a
basis for $H_1(\partial(N))$.

Let
$$
M=N\cup_h B^2 \times S^2
$$
be a Dehn filling of $N$ where $h:\partial(B^2 \times S^1) \to
\partial(N)$ is a homeomorphism such that $h(\partial(B^2) \times 0)$
is homologous to $\mu + n\lambda$ for some $n \in {\Bbb Z}$.

Now $h^{-1}(\partial(F) \times \{0,1/2\})$ bounds an annulus $A
\subset B^2 \times S^1$ which splits $B^2 \times S^1$ into two solid
tori $U_1$ and $U_2$ with $(U_i,A)$ homeomorphic to $(I \times I \times
S^1, I \times 0 \times S^1)$.

Then 
$$
V_1 = F \times [0,1/2]\cup_h U_1 \text{ and } V_2 = F \times
[1/2,1]\cup_h U_2
$$
are handlebodies of genus $2g(F)$, and thus give a splitting of $M$ which we call a
{\it horizontal splitting}.

If $\phi$ is (isotopic to) a periodic homeomorphism, then $N$ is
Seifert fibered, and if $n \ne 0$ this extends to a Seifert fibration
of $M$ -- with a singular fiber in $M-N$ if $n \ne \pm 1$. The other
singular fibers correspond to fixed points of some power  $\phi^k$ of
$\phi$ where $k$ is some proper divisor of the period $p$ of
$\phi$. The base surface of this Seifert fibration will be $\hat
F/\phi$ where $\hat F$ is obtained from $F$ by capping off
$\partial(F)$ with a 2-cell (on which $\phi = 1$).

Using the Riemann--Hurwitz formula and 2.1 one can show that the genus
of this horizontal splitting  is less than the genus of any vertical splitting
only when  $\hat F / \phi = S^2$ and $p=2$ or there is at most one
singular fiber (and $\hat F \ne S^2$). Moreover,  the only additional
cases with ``horizontal genus'' = `` vertical genus'' occur with $F^ =
T^2$. With a bit more analysis, including the classification of
periodic homeomorphisms of $T^2$ [H1;Theorem 12.11], one can establish

\proclaim{4.1 Theorem} Let $M$ be a closed oriented Seifert fibered
space over an oriented surface $S$ which has a horizontal splitting of
genus $g_{hor}$. Let $g_{ver}$ be the minimal genus of a vertical
splitting of $M$. Then

If $g_{hor} < g_{ver}$, then  either

(1.1) $M$ is Seifert fibered over $S^2$ with an even number $m \ge 4$
      of singular fibers, Seifert invariants $1/2,1/2,\dots,
      1/2,n/(2n+1)\quad; n\ge 1$, and Euler number $e = m/2$, or

(1.2) $M$ is Seifert fibered over a surface of genus $g > 0$ with at
      most one singular fiber and non-normalized invariant $\pm 1/n;
      \quad n \ge 1$.

If $g_{hor}= g_{ver}$, then $M$ is Seifert fibered over $S^2$ with
    $m=3$ singular fibers, Euler number $e = 1$ and invariants either

(2.1) $1/2,1/3,n/(6n\pm 1);\quad n \ge 0$, or

(2.2) $1/2,1/4, n/4n \pm 1); \quad n \ge 0$, or

(2.3) $1/3,1/3,n/(3n \pm 1); \quad n \ge 0$.
\endproclaim

In most of these cases we can show that the horizontal splitting of
$M$ is not represented by a positive diagram:

\proclaim{4.2 Theorem} Let $M$ satisfy the conclusion of 4.1. Then the
minimal genus horizontal splitting of $M$ is not represented by
a positive diagram in the cases:

(1.1) provided $n \ge 2$ or $m=4$,

(1.2) all cases, and

(2) all cases but $(1/2,1/3,1/5), (1/2,1/4,1/3),(1/3,1/3,1/2),\text
    { and } (1/2,1/3,1/7)$.
\endproclaim

\demo{Proof} First, the easiest case: (1.2). Here the horizontal
    splitting has genus $2g$ which is also the rank of $H_1(M)$. A
    positive diagram for this splitting would give a positive
    presentation for $\pi_1(M)$ with $2g$ generators and $2g$
    relations and thus a $2g \times 2g$ presentation matrix for
    ${\Bbb Z}^{2g}$ with all non-negative, and some positive
    entries. This is impossible.

Next we consider case (1.1). We note that identifying opposite edges on a
regular $2k$-gon $P$  (reversing orientation) produces an orientable
surface $\hat F$ of genus $g = k/2$ if $k$ is even or $g = (k-1)/2$ if
$k$ is odd.
Rotation of $P$ by $180\deg$ induces an orientation preserving
involution $\phi: \hat F \to \hat F$. Note that the fixed points of
$\phi$ come from the mid points of the edges of $P$, the center $c$ of
$P$, and, in case $k$ is even, the vertices of $P$. So the number of
singular fibers is $m = 2g+2$.

Let $F$ be obtained from $\hat F$ by removing an invariant
neighborhood of $c$. Then the manifold $M$ of (1.1) is obtained by
Dehn filling on the bundle $F \times [0,1]/(x,o) \sim (\phi(x),1)$.
Specifically, let $\lambda = \partial F \times 0$, as oriented by $F$ and
let $\mu$ be a transversal to $\lambda$ which traverses once in the
positive $t$ direction while traversing one-half turn in the positive
$\lambda$ direction (from $x_0 \in \partial F$ to $\phi(x_0)$). Then $M$
is obtained by the filling corresponding to $\mu\lambda^n = 1; \quad n
\in {\Bbb Z} - 0$.

Note that reflection through a diameter of $P$ induces an involution
of $F \times [0,1]/\phi$ which takes $\lambda$ to $\lambda^{-1}$ and $\mu
$ to $\mu\lambda^{-1}$. Thus  the surgeries $\mu\lambda^n =1$ and
$\mu\lambda^{-n-1} =1$ produce the same manifold (in case the reader
wonders if I have forgotten half of them). So we assume that $n \ge 1$.

The horizontal splitting gives a representation
$$
M = F \times [0,1] \cup_g F \times [0,1]
$$
where 
$$
g: \partial(F \times [0,1]) \to \partial(F \times [0,1])
$$
is a level preserving homeomorphism such that $g(x,0)=(x,0),
g(x,1)=(\phi(x),1)$ for all $x \in F$, and $g|\partial(F) \times [0,1]$
is a $n+1/2$ twist. We get a diagram $(S;X,Y)$ for this splitting as
follows. We put $S=\partial(F \times [0,1])$. We take arcs $a_1, \dots
,a_{2g}$ which cut $F$ to a 2-cell as shown in Figure 4 . Then $X =
\bigcup_i\partial (a_i \times [0,1])$ and $Y =\bigcup_i g(\partial (a_i
\times [0,1]))$.

\bigpagebreak
\epsfysize=2.5in
\centerline{\epsffile{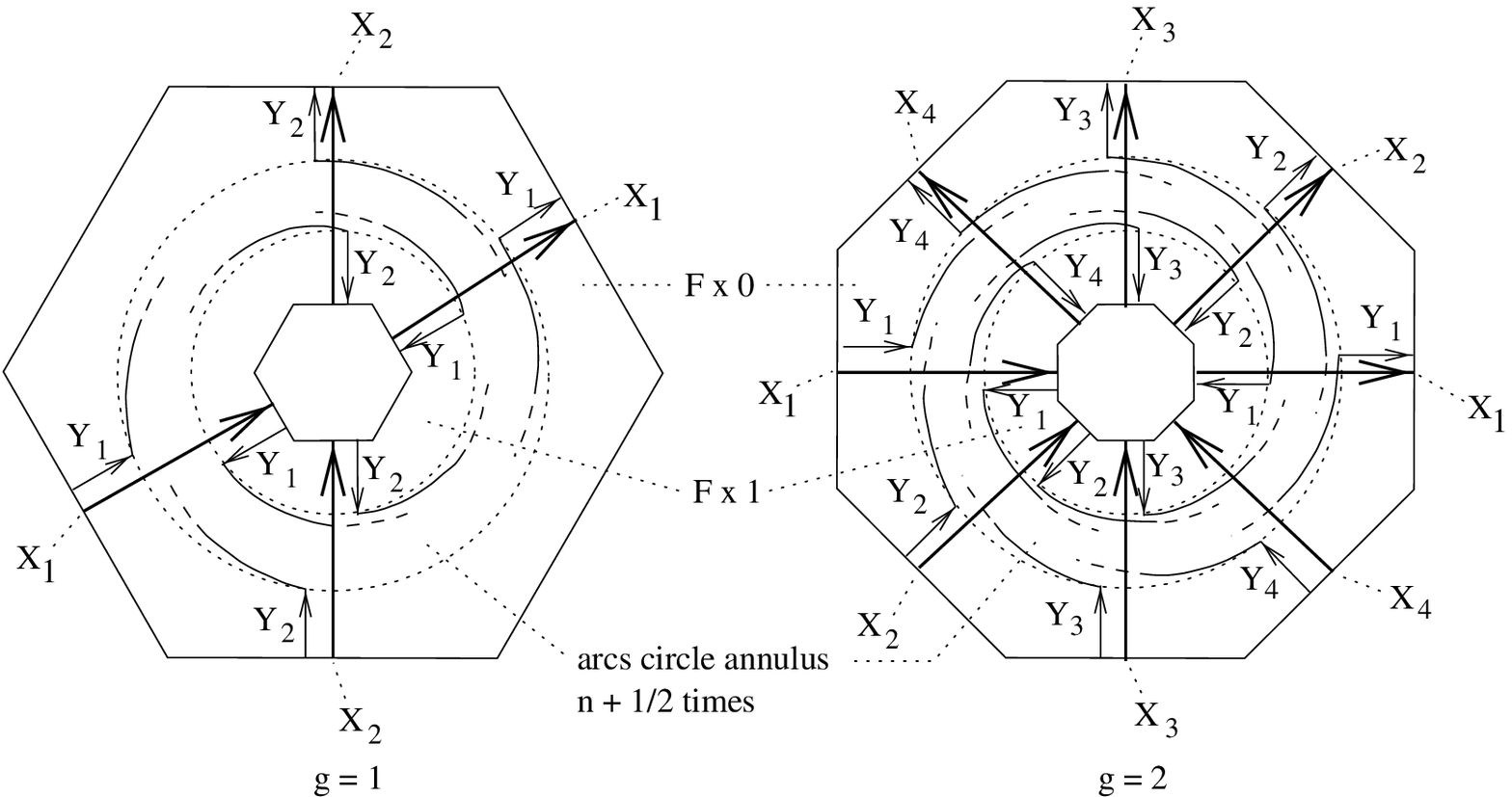}}
\botcaption{Figure 4}
Horizontal splitting for $(g=0; 1/2,1/2,\dots,1/2,n/2n+1; e=m/2)$: $1/n$
surgery on the
$S_g$ bundle, $g = m/2-1$, over $S^1$ with period two monodromy.
\endcaption
\bigpagebreak

Now most of ther components of $S-(X \cup Y)$ are squares with two
opposite sides in $X$ and the other two sides in $Y$. We form {\it
X-stacks} by taking maximal unions of such squares along common
$X$-edges, and similarly we form {\it Y-stacks}. See [H2] for more
details about stacks and their use in analysing Heegaard splittings. 

Now suppose $(S;X^*,Y^*)$ is a positive diagram for the splitting --
with $X^*$ meeting $X$ minimally, etc. Note that $(S;X,Y)$ has no waves. As described in [H2] each
component of $X^*$ must contain at least two $Y$-stack crossings and
each component of $Y^*$ must contain at least two $X$-stack crossings.

Now suppose that $n \ge 2$. Then if $a$ is any $X$-stack crossing and
$b$ any $Y$-stack crossing, $a\cup b$ separates $S$ into two
components each of which is incompressible in $F \times [0,1]$. Then
any component of $X^*$ other than the one containing $b$ must cross
$a$ twice in opposite directions. Thus $(S;X^*,Y^*)$ is not positive. 

Now suppose that $n=1$ and $m=4$. Then $g=1$. See Figure 5. Let $V_X$
(respectively $V_Y$) denote the handlebodies of the splitting with $X$
(respectively $Y$) as meridians. Suppose $D \subset V_X$ and $E
\subset V_Y$ are properly embedded, essential, oriented disks such
that $<\partial(D),\partial(E)>_p = +1$ at every point $p \in \partial(D) \cap
\partial(E)$. We may suppose that $X$ meets $\partial(D)$ and $Y$ meets
$\partial(E)$ minimally.

\bigpagebreak
\epsfysize=2.3in
\centerline{\epsffile{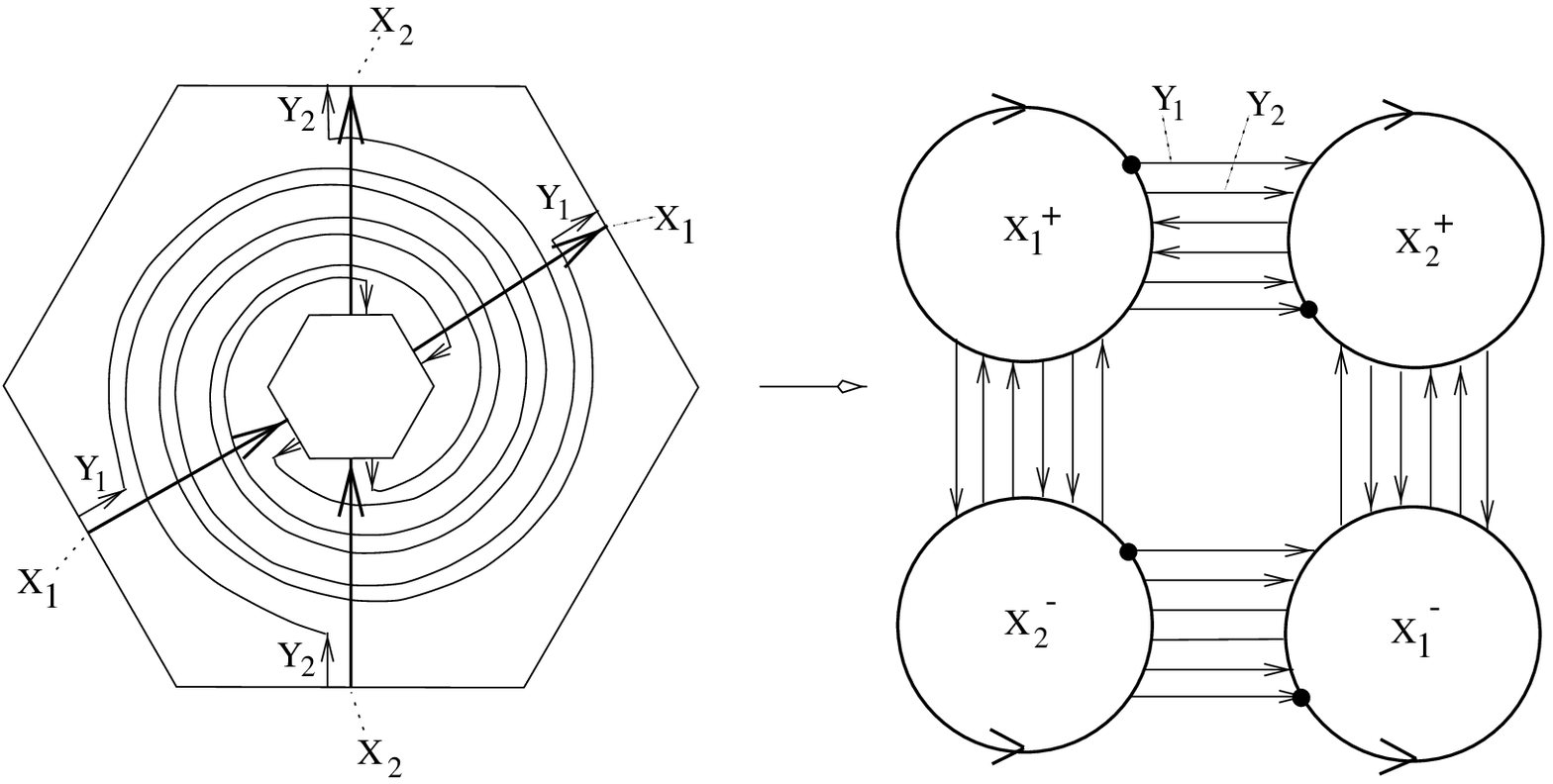}}
\botcaption{Figure 5}
Horizontal splitting for $(g=0;1/2,1/2,1/2,1/3;e=2)$.
\endcaption
\bigpagebreak

\proclaim{Claim} $\partial(D)$ ($\partial(E)$) contains crossings of two
distinct $Y$-stacks ($X$-stacks).
\endproclaim

We complete the proof of this case modulo the claim -- whose proof is
then given. So let $a_1,a_2 \subset \partial(D)$ ($b_1,b_2 \subset
\partial(E)$) be the $Y$-stack ($X$-stack) crossings from the
claim. Because they are crossings of different stacks, we see that
$J = a_1\cup a_2 \cup b_1 \cup b_2$ separates $S$ into two incompressible
components.

Now suppose that $D$ and $E$ come from a positive diagram for the
given splitting. There is a second disk $E' \subset V_Y$ comming from
the diagram. Each $a_i$ crosses $b_1$ and $b_2$ in the same direction;
so $a_1$ and $a_2$ are consistantly oriented in $J$. Now $\partial(E')$ must cross $J$ at least two times in
opposite directions. Since $E'\cap (b_1 \cup b_2)= \emptyset $,
$\partial(E')$ must cross $a_1 \cup a_2$ at least two times in
opposite directions. So the diagram was not positive after all.

{\bf Proof of Claim.} The arguement is symmetric in $X$ and $Y$. Let
$D_1,D_2 \subset V_X$ be the disks bounded by the components of
$X$. Cutting $V_X$ along $D_1 \cup D_2$ produces a 3-cell $B$ whose
boundary contains two copies $D_{i^+}$ and $D_{i^-}$ of $D_i; \quad
i=1,2$. We get a graph $\Gamma \subset \partial(B)$ whose ``fat''
vertices are these $D_{i^\pm}$ and whose edges are the arcs of
$\partial(D) $ cut open.

Now $\Gamma$ has at least two loops -- comming from outermost
components of $D-(D_1 \cup D_2)$. Suppose that one of these loops , $a$ is
based at, say,  $D_{1^+}$. This loop must separate one ({\it the
singleton}) of the remaining vertices from the other two. If the singleton is $D_{1^-}$,
then $a$ must cross both of the $Y$-stacks with one side in $D_{1^-}$,
and we are done.

So, say, the singleton is $D_{2^+}$. Let $n_{i,j}; \quad i,j \in
\{1^\pm,2^\pm\}$ denote the number of edges of $\Gamma$ with one
vertex in $D_i$ and the other in $D_j$. So $n_{1^+,1^+} >0$, and
$n_{1^-,2^+} = n_{2^-,2^+}=n_{2^+,2^+} =0$.

The edges of $\Gamma$ must match up on reidentification. This forces
consistency equations:
$$
n_{1^-,1^-}+n_{1^-,2^-}+n_{1^-,2^+}
=n_{1^+,1^+}+n_{1^+,2^-}+n_{1^+,2^+},
$$
$$
n_{1^-,2^-}+n_{1^+,2^-}+n_{2^-,2^-}=n_{1^-,2^+}+n_{1^+,2^+}+n_{2^+,2^+}
$$
If $n_{1^-,1^-}=0$, then putting in all of the $0$'s and solving for
$n_{1^-,2^-}$ from the first equation and putting this  in the second
gives
$$
n_{1^+,1^+}+2n_{1^+,2^-}+n_{2^-,2^-} =0
$$
This is impossible, as all terms are non-negative and $n_{1^+,1^+}
>0$. So $n_{1^-,1^-} > 0$. This forces the singleton to be $D_{1^-}$
and  completes the proof of the claim.

Finally we consider the case (2). Here $M$ is a torus bundle over
$S^1$ (see [H1; Theorem 12.11]), and is obtained as above with $\phi$ the
rotation of the hexagonal torus by $60^o$ or $120^o$ or rotation
of the square torus by $90^o$ (the period $2$ rotation of either
gives rise to four singular fibers and is included in case (1.2)
above). Figure 6  shows diagrams for the corresponding splittings. The
right hand side shows the diagrams in a more standard form and are
valid for $n >0$. With
notation  as in (1.2) the manifolds come from $\mu\lambda^k=1$ Dehn
filling on $F \times [0,1]/\phi$ for some $k \in {\Bbb Z} -0$. Here we do not have duality between
positive and negative $k$. The $n$ in the statement is $n = |k|$ and
the sign in the denominator of the third invariant is $+$ or $-$
according as $k >0$ or $k<0$.

\bigpagebreak

\topcaption{Figure 6}
Horizontal splittings for the $1/\pm n$ surgeries on the torus bundles over
$S^1$ with periodic monodromy ($\ne 2$).
\endcaption

\epsfysize=2.3in
\centerline{\epsffile{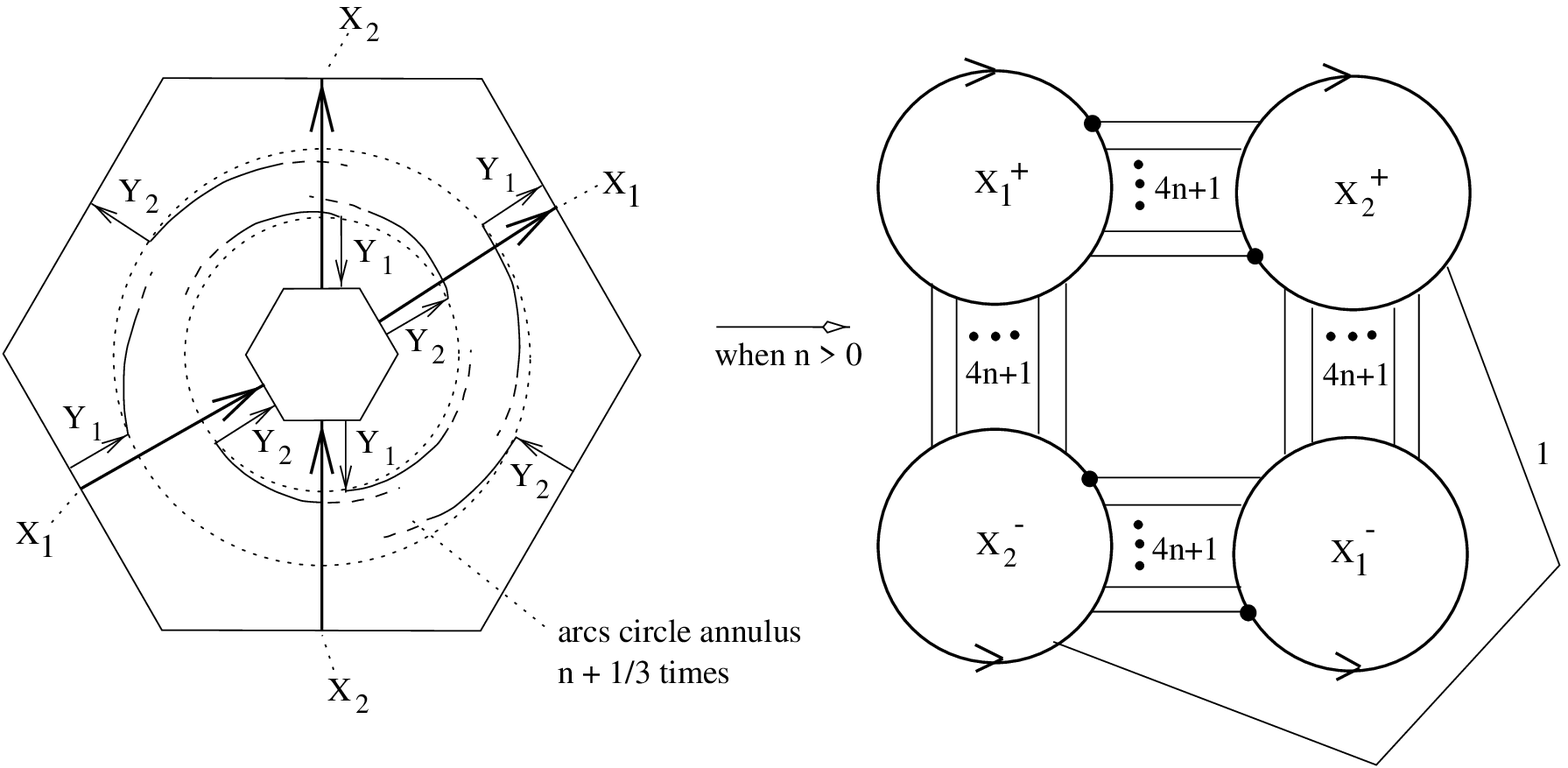}}
\botcaption{Figure 6a}
Period three: $(g=0,(1/3,1/3,n/(3n\pm1)); e=1)$.
\endcaption
\bigpagebreak

\bigpagebreak
\epsfysize=2.3in
\centerline{\epsffile{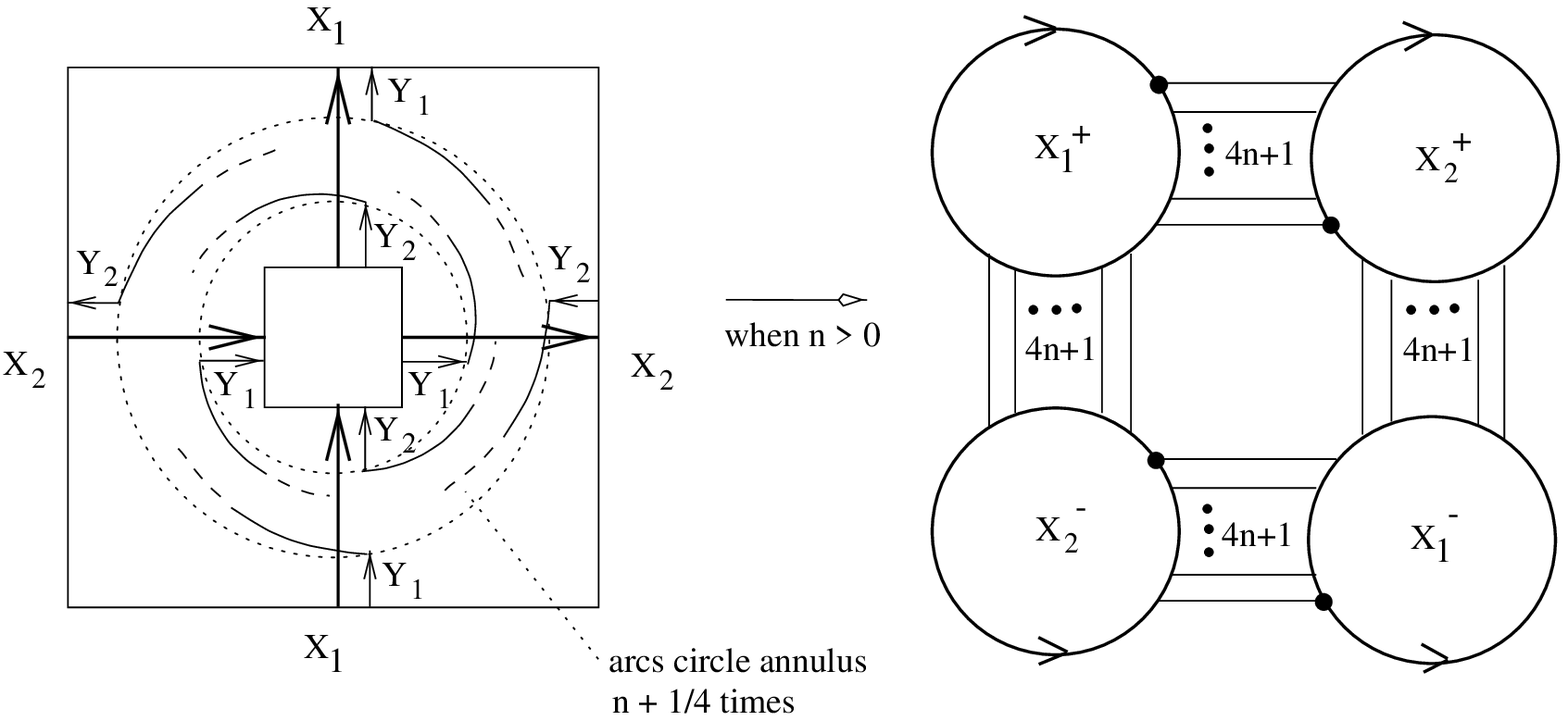}}
\botcaption{Figure 6b}
Period four: $(g=0,(1/2,1/4,n/(4n\pm1)); e=1)$.
\endcaption
\bigpagebreak

\bigpagebreak
\epsfysize=2.3in
\centerline{\epsffile{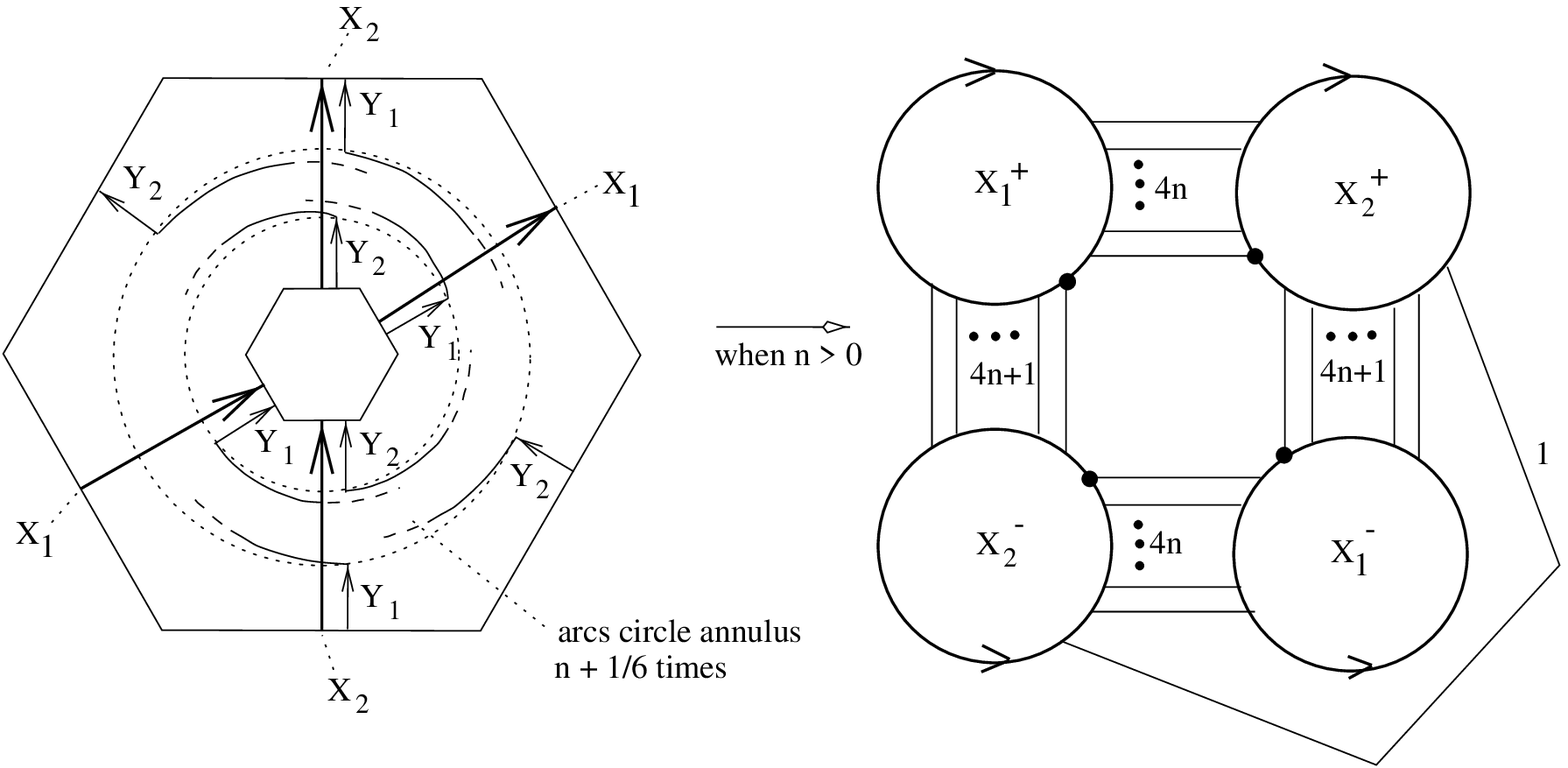}}
\botcaption{Figure 6c}
Period three: $(g=0,(1/2,1/3,n/(6n\pm1)); e=1)$.
\endcaption
\bigpagebreak

The cases $k = 0$ all give lens spaces. In the cases $k = -1$, the
horizontal splitting is also a vertical splitting (by, for example,
[H2; Theorem 3.2]) and so is
represented by a positive diagram and is excluded. The remaining
possibilities follow very much as in the proof of  case (1.2); with
$k = 1,- 2$, and $k=2$ when the period $p$ of $\phi $ is six
requiring the extra arguements about crossings of  two distinct stacks. This works except
when $p=6$ and $k=1$. Here not every
$X$-stack meets every $Y$-stack; so we cannot claim the properties of
$J = a_1\cup a_2\cup b_1\cup b_2$. This accounts for the final
exclusion. \qed
\enddemo

{\bf 5. Remarks.} There are some questions left open. With respect to
Seifert manifolds (orientable with orientable base):

(1) We do not know whether the horizontal splittings of the manifolds
$$
g=0, m \ge 6 \text{ and even, invariants } = 1/2, \dots, 1/2, 1/3, e =
m/2
$$
or
$$
g=0, m=3,\text{ invariants } = 1/2/1/3/1/7, e=1
$$
are represented by positive diagrams.

(2) We do not know whether the minimal genus vertical splittings with
$$
g > 0 \text{ and } m \le 2
$$
are represented by positive diagrams.

More generally

(3) We do not know whether all vertical splittings ( when $g >0$) are
    represented by positive diagrams.

Outside the class of Seifert manifolds we know nothing.

%%%%%%%%%%%%%%%%%

\enddocument
\end